\documentclass[10pt]{amsart}
\usepackage{amssymb,amsbsy,amsmath,amsfonts,amssymb}
\usepackage{latexsym,euscript,exscale}

\title{Infinite asymptotic games}
\author {Christian Rosendal}
\date {April 2006}
\newcommand{\ram}{[\N]^{\N}}

\newcommand{\ra}[1]{[#1]^{\N}}

\newcommand {\A}{\mathbb A}
\newcommand {\B}{\mathbb B}

\newcommand {\N}{\mathbb N}

\renewcommand{\leq}{\ensuremath{\leqslant}}
\renewcommand{\geq}{\ensuremath{\geqslant}}

\newcommand{\PP}{\mathbb P}

\newcommand {\E}{\mathbb E}

\newcommand{\eps}{\epsilon}

\newcommand{\iso}{\cong}

\newcommand{\tom} {\emptyset}

\newcommand{\hviss}{\leftrightarrow}

\newcommand{\saa}{\Rightarrow}
\newcommand{\equi}{\Leftrightarrow}

\newcommand{\til}{\rightarrow}

\newcommand{\Lim}[1]{\mathop{\longrightarrow}\limits_{#1}}

\newcommand {\del}{ \; \big| \;}

\newcommand {\ku} {\mathcal}

\newcommand {\e} {\exists}
\renewcommand {\a} {\forall}

\newcommand{\fed}{\boldsymbol}

\newtheorem{thm}{Theorem}
\newtheorem{cor}[thm]{Corollary}
\newtheorem{lemme}[thm]{Lemma}
\newtheorem{prop} [thm] {Proposition}
\newtheorem{defi} [thm] {Definition}

\begin{document}
\thanks{The  author is partially supported by NSF grant DMS 0556368}

\subjclass[2000]{Primary: 46B03, Secondary 03E15}

\keywords{Infinite asymptotic games, extraction of subsequences,
weakly null trees.}

 \maketitle

\begin{abstract}
We study infinite asymptotic games in Banach spaces with an F.D.D.
and prove that analytic games are determined by characterising
precisely the conditions for the players to have winning strategies.
These results are applied to characterise  spaces  embeddable into
$\ell_p$ sums of finite dimensional spaces, extending results of
Odell and Schlumprecht, and  to study various notions of homogeneity
of bases and Banach spaces. The results are related to questions
of rapidity of subsequence extraction from normalised weakly null
sequences.
\end{abstract}
\tableofcontents
\section{Introduction}
A number of results have surfaced over some years that involve questions about
Banach spaces of the following kind: Suppose $E$ is a Banach space such that
every normalised weakly null sequence has a subsequence with a certain
property. What can then be concluded about $E$? In general, it is not enough
that one can just find some subsequence, but in various guises (for example,
regulators \cite{fonf}, weakly null trees \cite{kalton,os}) it has been noticed
that if the subsequence can be chosen sufficiently fast then this is sufficient
to provide some information about the whole space.

The notion behind this principle was crystalised through work of E. Odell and
T. Schlumprecht, in particular \cite{os} (or see \cite {os2} for a survey), and
can be formulated via {\em infinite asymptotic games}. Infinite asymptotic
games, as defined in, e.g., \cite{os} and \cite{fr}, are a straightforward
generalisation of the finite asymptotic games of B. Maurey, V. Milman, and N.
Tomczak-Jaegermann in \cite{mmt}, but with very different objectives. While
finite asymptotic games are used for describing the asymptotic
finite-dimensional structure of Banach spaces, infinite asymptotic games have a
much more combinatorial flavour. For example, they can be used as basis for a
simplified proof of Gowers' block Ramsey principle as in \cite{exact} and are
also related to the notion of {\em tight} Banach spaces from \cite{tight}.

A {\em weakly null tree} in a Banach space can be defined as a normalised
sequence $(x_s)_{s\in \N^{<\infty}}$ such that for each $s$, the sequence
$(x_{sn})_{n\in \N}$ is weakly null (we shall use a trivial variation of this
definition later on). Odell and Schlumprecht proved that any reflexive
separable Banach space in which any weakly null tree has a branch
$C$-equivalent to $\ell_p$ embeds into an $\ell_p$-F.D.D. Their proof involved
first proving the equivalence of their hypothesis with the existence of a
winning strategy for I in a {\em closed} infinite asymptotic game.

The main result of this paper, Corollary \ref{strategy for I}, is to give an
exact criterion for when I has a winning strategy in any {\em coanalytic}
infinite asymptotic game, thus extending the result of Odell and Schlumprecht.
The proof of Corollary \ref{strategy for I} is very different from the proof of
the case of closed games and is instead related to a result in \cite{fr}.

The game behind Gowers' block Ramsey principle mentioned above is a a game for
extracting block sequences of block sequences (see Bagaria and L\'opez-Abad
\cite{jordi}). Analogously, we introduce games for extracting subsequences of
block sequences and show that the strongest of these games is actually
equivalent with the infinite asymptotic game, Theorem \ref{equivalence of
games}. We also study weaker versions of these games in Section
\ref{subsequence extraction} indicating their differences with the infinite
asymptotic game.

Corollary \ref{strategy for I} allows us to use the techniques of Odell and
Schlumprecht in a priori more complicated settings. For example, we prove that
if $E$ is a reflexive separable Banach space in which any weakly null tree has
a subsymmetric branch, then $E$ embeds into some $\ell_p$-F.D.D., Theorem
\ref{os}. Also, assuming strong axioms of set theory, we show that for a space
$E$ which is not $\ell_1$-saturated, every weakly null tree in $E$ has a branch
spanning a space isomorphic to $E$ if and only if $E$ has an unconditional
block homogeneous basis, i.e., an unconditional basis all of whose blocks span
isomorphic spaces, Theorem \ref{block homo}.

For the reader not familiar with descriptive set theory,
we have included subsection \ref{descriptive} for a brief review of used notions and results.
\section{Theory}\label{theory}

\subsection{Infinite asymptotic and subsequence games in vector spaces}
Let $E$ be a real or complex vector space with basis $(e_i)_{i\in \N}$. For a
vector $x=\sum_ia_ie_i\in E$, we designate by ${\rm supp}(x)$ the finite set
$\{i\in \N\del a_i\neq 0\}$ and use $k< x< n$ to denote that $k<{\rm supp}(x)<
n$, i.e., that $k< \min{\rm supp}(x)$ and $\max{\rm supp}(x)< n$. Similar
notation is used for inequalities between non-zero vectors. A finite or
infinite {\em block sequence} is a sequence $(x_i)$ of non-zero vectors such
that $x_i<x_{i+1}$ for all $i$. We denote the set of infinite, resp. finite,
block sequences of $E$ by $E^\infty$, resp. $E^{<\infty}$. If $(x_i)$ is a
finite or infinite sequence of vectors, we denote by $[x_i]$ the corresponding
linear span, On the other hand, when $(x_i)$ is a sequence in a {\em Banach}
space, we let $[x_i]$ be the closed linear span of the vectors.

The {\em infinite asymptotic game} (IAG) on $E$ between two players I and II is
defined as follows: I and II alternate (with I beginning) in choosing
respectively natural numbers $n_0, n_1, n_2,\ldots$ and non-zero vectors
$x_0<x_1<x_2<\ldots\in E$ according to the constraint $n_i< x_i$:
$$
\begin{array}{cccccccccccc}
{\bf I} & & & n_0 &  & n_1 &  & n_2 & & n_3& &\ldots \\
{\bf II} & & &  & n_0<x_0 &  & n_1<x_1 &  & n_2<x_2 & & n_3<x_3 &\ldots
\end{array}
$$
We say that the sequence $(x_n)_{n\in \N}$ is the {\em outcome} of the game.

There is another natural game to play on a basis. This is the {\em subsequence
game} (SG) defined as follows. Player I and II alternate in choosing
respectively digits $\eps_0,\eps_1,\eps_2,\ldots\in\{0,1\}$ and non-zero
vectors $x_0<x_1<x_2<\ldots\in E$, now with II beginning:
$$
\begin{array}{cccccccccccc}
{\bf I} & & &  & \eps_0 &  & \eps_1 &  & \eps_2 & & \eps_3 &\ldots\\
{\bf II} & & & x_0 &  & x_1 &  & x_2 & & x_3& &\ldots
\end{array}
$$
We thus see II as constructing an infinite block sequence $(x_n)_{n\in \N}$,
while I chooses a subsequence $(x_n)_{n\in A}$ by letting $n\in A\equi
\eps_n=1$. This subsequence $(x_n)_{n\in A}$ is then called the {\em outcome}
of the game.

Despite their superficial difference, we shall prove that the games are
actually {\em equivalent}, i.e., for any $\A\subseteq E^\infty$, player I,
resp. II, has a strategy in (IAG) to play in $\A$ if and only if I, resp. II,
has a strategy in (SG) to play in $\A$.

A {\em block tree} is a non-empty infinitely branching tree $T\subseteq
E^{<\infty}$ such that for all $(x_0,x_1,\ldots,x_n)\in T$ the set
$$
\{y\in E\del (x_0,x_1,\ldots,x_n,y)\in T\}
$$
can be written as $\{y_j\}_{j\in \N}$ for some infinite block sequence
$(y_j)_{j\in \N}\in E^\infty$ such that $x_n<y_0$. We denote by $[T]$ the set
of {\em infinite branches} of $T$, i.e., the set of $(x_i)\in E^\infty$ such
that $(x_0,\ldots,x_n)\in T$ for all $n$.

\begin{thm}\label{equivalence of games}
The games (IAG) and (SG) are equivalent. Moreover, for any $\A\subseteq
E^\infty$, II has a strategy in (IAG) to play in $\A$ if and only if there is a
block tree $T$ such that $[T]\subseteq \A$.
\end{thm}

\begin{proof}Let $\A\subseteq E^\infty$ be given.
Suppose first that I has a strategy $\sigma$ in (IAG) to play in $\A$. We
construct a strategy for I in (SG) to play in $\A$. So suppose II is playing a
sequence $x_0<x_1<x_2<\ldots$ in (SG). For every $k\geq 0$, we have to
construct an answer $\eps_k$ by I based only on the first $k+1$ vectors
$x_0<\ldots<x_k$. So let $n_0$ be the first number played by I according to
$\sigma$ in (IAG). I then plays $\eps_0=\eps_1=\ldots=0$ until $k_0$ is minimal
such that $n_0<x_{k_0}$ and then lets $\eps_{k_0}=1$. Let now $n_1$ be the
answer of $\sigma$ to $x_{k_0}$ played by II. I then plays
$\eps_{k_0+1}=\eps_{k_0+2}=\ldots=0$ until $k_1>k_0$ is minimal such that
$n_1<x_{k_1}$, at which point $\eps_{k_1}$ is set to be $1$. Again we let $n_2$
be the answer by $\sigma$ to $(x_{k_0},x_{k_1})$ played by II in (IAG) and I
continues as before. Since $\sigma$ forces the outcome $(x_{k_i})_{i\in \N}$ to
be in $\A$, this describes a strategy in (SG) for I to play in $\A$.

Suppose conversely that I has a strategy $\sigma$ in (SG) to play in $\A$,
i.e., $\sigma$ associates to each finite block sequence $(x_0,\ldots,x_n)$ a
digit $\eps_n\in\{0,1\}$. We define a tree of {\em good} finite block sequences
and a function $\phi$ associating to each good $(x_0,\ldots,x_n)$ another block
sequence $(y_0,\ldots,y_m)$, $n\leq m$, such that $(x_0,\ldots,x_n)$ is the
subsequence of $(y_0,\ldots,y_m)$ extracted by $\sigma$. We begin by letting
$\tom$ be good with $\phi(\tom)=\tom$. Now if $(x_0,\ldots,x_n)$ is good with
$\phi(x_0,\ldots,x_n)=(y_0,\ldots,y_m)$, we let $(x_0,\ldots,x_n,x_{n+1})$ be
good if there is $(y_{m+1},\ldots,y_k)$ such that
$$
\sigma(y_0,\ldots,y_m,y_{m+1},\ldots,y_l)=\left\{
                                            \begin{array}{ll}
                                              0, & \hbox{if $m<l<k$;} \\
                                              1, & \hbox{if $l=k$.}
                                            \end{array}
                                          \right.
$$
Moreover, we let
$\phi(x_0,\ldots,x_n,x_{n+1})=(y_0,\ldots,y_m,y_{m+1},\ldots,y_k)$ for some
such $(y_{m+1},\ldots,y_k)$. Now if $(x_0,x_1,x_2,\ldots)$ is such that
$(x_0,\ldots,x_n)$ is good for all $n$, then $\phi(\tom)\subseteq
\phi(x_0)\subseteq \phi(x_0,x_1)\subseteq \ldots$ and $(x_i)$ is the
subsequence extracted by $\sigma$ from $\bigcup_m\phi(x_0,x_1,\ldots,x_m)$ and
hence belongs to $\A$. We claim that if $(x_0,\ldots,x_n)$ is good, then for
some $k\in \N$ if $k<y$, then $(x_0,\ldots,x_n,y)$ is good. Suppose not and let
$\phi(x_0,\ldots,x_n)=(y_0,\ldots,y_m)$. Then, we can find
$x_n=y_m<y_{m+1}<y_{m+2}<\ldots$ such that no $(x_0,\ldots,x_n,y_i)$ is good
and, in particular, $\sigma(y_0,\ldots,y_m,y_{m+1},\ldots,y_{k})=0$ for all
$k>m$. But then $(x_0,\ldots,x_n)$ is the subsequence of $(y_0,y_1,\ldots)$
extracted by $\sigma$, contradicting that $\sigma$ is a strategy for I to play
in $\A\subseteq E^\infty$. It is now trivial to use this claim to construct a
strategy for I in (IAG) to play in $\A$.

Similarly, if $\sigma$ is a strategy for II in (SG) to play in $\A$, we can
construct a strategy for II in (IAG) to play in $\A$ as follows. If I plays
$n_0$ in (IAG), we let I play $\eps_0=\eps_1=\ldots=0$ in (SG) and II respond
according to $\sigma$ until II plays some $x_{k_0}$ such that $n_0<x_{k_0}$. II
then plays $x_{k_0}$ as response to $n_0$ in (IAG) and I lets $\eps_{k_0}=1$.
Let $n_1$ be the next number played by I in (IAG) and let I play
$\eps_{k_0+1}=\eps_{k_0+2}=\ldots=0$ and II respond further using $\sigma$
until II plays some $x_{k_1}>n_1$, whence I responds by $\eps_{k_1}=1$ and II
plays $x_{k_1}$ in (IAG). We continue in this way and see that II plays
according to $\sigma$ in (SG) with I choosing the subsequence $(x_{k_i})$,
whereby $(x_{k_i})\in \A$. Moreover, $(x_{k_i})$ is exactly the sequence played
by II in (IAG) and hence II has a strategy to play in $\A$ in the game (IAG).

Suppose $\sigma$ is a strategy for II in (IAG) to play in $\A$,
i.e., $\sigma$ is a function associating to each finite non-empty
sequence $(n_0,\ldots, n_k)$ of natural numbers some $x_k>n_k$. We
define a pruned tree $T'$ of {\em good} finite block bases
$(x_0,\ldots,x_k)$ and a function $\psi$ associating to each good
$(x_0,\ldots,x_k)$ a sequence $(n_0,\ldots,n_k)$ such that
$(n_0,x_0,\ldots,n_k,x_k)$ is consistent with $\sigma$, i.e., for
all $l\leq k$, $\sigma(n_0,\ldots,n_l)=x_l$.
\begin{itemize}
  \item The empty sequence $\tom$ is good and $\psi(\tom)=\tom$.
  \item If $(x_0,\ldots,x_k)$ is good and
  $$
  \psi(x_0,\ldots,x_k)=(n_0,\ldots,n_k),
  $$
  then we let
  $(x_0,\ldots,x_k,y)$ be good if there is some $m>n_k$ such that
  $y=\sigma(n_0,\ldots,n_k,m)$ and in this case we let
  $$
  \psi(x_0,\ldots,x_k,y)=(n_0,\ldots,n_k,m'),
  $$
  where $m'$ is the least such
  $m$.
\end{itemize}
Now, if $(x_0,x_1,x_2,\ldots)$ is such that $(x_0,\ldots,x_k)$ is good for all
$k$, then $\psi(\tom)\subseteq \psi(x_0)\subseteq \psi(x_0,x_1)\subseteq
\ldots$ and $x_k=\sigma(\psi(x_0,\ldots,x_k))$ for all $k$, whence
$(x_0,x_1,x_2,\ldots)$ is a play of II according to $\sigma$ and hence belongs
to $\A$. So $[T']\subseteq \A$. Also by construction, for each
$(x_0,\ldots,x_k)\in T'$ and $n$ there is $y>n$ such that
$(x_0,\ldots,x_k,y)\in T'$ and thus it is easy to construct a block subtree
$T\subseteq T'$, whereby also $[T]\subseteq \A$.

If $T$ is a block tree all of whose branches lie in $\A$, then we can construct
a strategy for II in (SG) to play in $\A$ as follows: First play
$x^{(0)}_0<x^{(1)}_0<\ldots$ such that $(x^{(i)}_0)\in T$ for all $i$ until
$k_0$ is minimal with $\eps_{k_0}=1$. Then play $x^{(0)}_1<x^{(1)}_1<\ldots$
such that $(x^{(k_0)}_0,x^{(i)}_1)\in T$ for all $i$ until $k_1\geq 0$ is
minimal with $\eps_{k_0+k_1}=1$, etc. Then the subsequence $(x^{(k_i)}_i)_{i\in
\N}$ chosen by I will be a branch of $T$ and hence lie in $\A$.
\end{proof}

The fact that strategies for II in (IAG) can be refined to a block tree is
perhaps not too surprising, as the outcome of the game is independent of the
play of I. However, there is a small twist, as, in general, taking $T'$ to be
the set of all $(x_0,\ldots,x_k)$ that are played according to $\sigma$ does
not produce a tree all of whose branches lie in $\A$. This is avoided above by
using ideas of D. A. Martin \cite{martin}.

\subsection{Notation and results from descriptive set
theory}\label{descriptive} We recall that a {\em Polish} space is a separable
topological space whose topology can be induced by a complete metric. The {\em
Borel} sets are those belonging to the smallest $\sigma$-algebra containing the
open sets and a {\em standard Borel} space is the measurable space obtained by
equipping the underlying set of a Polish space with its $\sigma$-algebra of
Borel sets. A subset $A\subseteq X$ of a Polish or standard Borel space is {\em
analytic} or ${\fed \Sigma}_1^1$ if there is a Polish space $Y$ and a Borel set
$B\subseteq X\times Y$ such that
$$
x\in A\equi \e y\in Y\; (x,y)\in B.
$$
In other words, $A$ is the projection ${\rm proj}_X(B)$ of $B$ onto $X$. All
Borel sets are analytic. A set $C\subseteq X$ is {\em coanalytic} or
${\fed\Pi}_1^1$ if its complement is analytic, or, equivalently, if there is a
Polish space $Y$ and a Borel set $D\subseteq X\times Y$ such that
$$
x\in C\equi \a y\in Y\; (x,y)\in D.
$$
And if $n\geq 1$, a set $A$ is ${\fed \Sigma}_{n+1}^1$ if there is a ${\fed
\Pi}_n^1$ set $B\subseteq X\times Y$ for some Polish space $Y$ such that
$A={\rm proj}_X(B)$. Also, $C$ is ${\fed \Pi}_n^1$ if its complement is ${\fed
\Sigma}_n^1$.  Sets belonging to some class ${\fed \Sigma}_n^1$ or ${\fed
\Pi}^1_n$ are called {\em projective} and encompass most of the sets used in
analysis. For example, if $E$ is a separable Banach space, then the set of
sequences in $C[0,1]$, whose closed linear span is isomorphic to $E$, is
analytic in $C[0,1]^\N$.

If $A$ is an infinite subset of $\N$, we denote by $[A]^\N$ the set of all
infinite increasing sequences in $A$ or equivalently the infinite subsets of
$A$. Also, if $n_1<\ldots<n_k$ are natural numbers, we denote by
$[(n_1,\ldots,n_k),A]$ the elements of $[A]^\N$ whose first $k$ terms are
$n_1,\ldots,n_k$. We give $[A]^\N$ the Polish topology inherited from its
inclusion in $A^\N$, where $A$ is taken discrete. A theorem due to F. Galvin
and K. Prikry says that if $\A\subseteq [\N]^\N$ is Borel, then there is an
infinite set $A\subseteq \N$ such that either $[A]^\N\cap \A=\tom$ or
$[A]^\N\subseteq \A$. This was improved by J. Silver to include the
$\sigma$-algebra $\sigma({\bf \Sigma}_1^1)$ generated by the analytic sets.

A result due to D.A. Martin states that all countable games on integers with
Borel winning conditions are determined. The set theoretical statement called
{\em Projective determinacy} says that all games on integers with a projective
winning condition is determined. This is not provable from the usual axioms of
set theory, but in contrast with other strong axioms of set theory is a part of
what many consider to be the right axioms of set theory. It has the consequence
of extending the regularity properties of ${\fed \Sigma}_1^1$ to all projective
sets. Thus, a result of L. Harrington and A.S. Kechris says that the
Galvin--Prikry result extends to all projective sets under projective
determinacy.  We refer to \cite{kechris} for these results.

We use $(x_n)\sim(y_n)$ to denote that two sequences in Banach spaces $X$ and
$Y$ are {\em equivalent}, i.e., that the mapping $T\colon x_n\mapsto y_n$
extends to an isomorphism of their closed linear spans. Also $X\iso Y$ denotes
that the spaces $X$ and $Y$ are isomorphic. If we index these relations by a
constant $K$, we mean that the constant of equivalence or isomorphism is
bounded by $K$. If $X$ is a Banach space, we denote its unit sphere by $\ku
S(X)$ or $\ku S_X$.

If $X$ is a Banach space and $(F_i)_{i\in \N}$ a sequence of finite-dimensional
subspaces, we say that $(F_i)$ is a {\em finite-dimensional decomposition} or
F.D.D. of $X$ if any element $x\in X$ can uniquely be written as a
norm-convergent series $x=\sum_if_i$, where $f_i\in F_i$.

\subsection{Gowers' block sequence game}
The infinite asymptotic game resembles W.T. Gowers' {\em block sequence game}
\cite{gowers} that is defined as follows: Player I and II alternate in choosing
infinite-dimensional subspaces $Z_0,Z_1,\ldots\subseteq E$ and non-zero vectors
$x_1<x_2<\ldots$ in $E$ according to the constraint $x_i\in Z_i$:
$$
\begin{array}{cccccccccccccc}
{\bf I} & & Z_0 & &Z_1&  & Z_2 &  & \ldots \\
{\bf II}   &  &  & x_0\in Z_0 &  & x_1\in Z_1 &&x_2\in Z_2 &\ldots
\end{array}
$$
Again, the block sequence $(x_i)$ is called the {\em outcome} of the game.

Gowers \cite{gowers} essentially proved that if $\A\subseteq E$ is an analytic
set such that any  infinite-dimensional $Z\subseteq E$ contains a block
sequence belonging to $E$, then there is an infinite-dimensional subspace
$X\subseteq E$ such that II has a strategy in the block sequence game, in which
I is restricted to playing subspaces of $X$, to play inside a slightly bigger
set than $\A$.

We shall prove a similar result, Theorem \ref{main}, for the infinite
asymptotic game, except that there is no passage to a subspace $X$ involved.
This implies that the nature of the infinite asymptotic game is really
different from that of Gowers' game. On the other hand, in \cite{exact} it is
shown that if $E$ is a vector space over a countable field, $\A\subseteq
E^\infty$ is analytic, and for all infinite-dimensional $X\subseteq E$, II has
a strategy in (IAG) played below $X$ to play in $\A$, then there is $X$ such
that II has a strategy in the block sequence game, in which I is restricted to
playing subspaces of $X$, to play in $\A$. This elucidates the exact relation
between the two games.

\subsection{Infinite asymptotic games in normed vector spaces}\label{evn}
Suppose now that $E$ is moreover a normed vector space and $(e_i)$ is a
normalised basis. Since each $[e_0,\ldots,e_n]$ is a finite-dimensional Banach
space, it is complete and hence Polish. We can therefore naturally see $E$ as a
standard Borel space by letting it be the increasing union of the sequence of
Borel subsets $[e_0,\ldots,e_n]\subseteq E$. We denote by $bb(e_i)$ and
$fbb(e_i)$ the set of infinite, resp. finite, {\em normalised} block sequences
and notice that the former is a Borel subset of $E^\N$.

When playing the infinite asymptotic or subsequence game in the normed space
$E$, we will now demand that II always plays normalised vectors. Also, block
trees are now supposed to consist of normalised vectors. It is obvious that the
proof of Theorem \ref{equivalence of games} adapts to this context.

 Let also $\A\subseteq bb(e_i)$
be non-empty and $\Delta=(\delta_i)_{i=0}^\infty$ be a decreasing sequence of
strictly positive reals converging to $0$ (which we denote simply by
$\Delta>0$). We let
$$
\A_\Delta=\{(y_i)\in bb(e_i)\del \e (x_i)\in \A\; \a i\;
\|x_i-y_i\|<\delta_i\},
$$
and
\begin{displaymath}\begin{split}
{\rm Int}_\Delta(\A)=&\{(y_i)\in bb(e_i)\del \a (x_i)\in bb(e_i)\;( \a i\; \|x_i-y_i\|<\delta_i\til (x_i)\in \A)\}\\
=&\sim(\sim\A)_\Delta.
\end{split}\end{displaymath}
Thus ${\rm Int}_\Delta(\A)\subseteq \A\subseteq \A_\Delta$. We should notice
here that if $(e_i)$ is a Schauder basis for the completion of $E$ and $\A$ is
a set which is closed under taking equivalent sequences, i.e., if for
$(x_i)\sim (y_i)$ we have $(x_i)\in \A\hviss(y_i)\in \A$, then for $\Delta>0$
chosen sufficiently small (depending on the constant of the basis $(e_i)$), we
have ${\rm Int}_\Delta(\A )=\A=\A_\Delta$. Finally, define the following
relation $\ku R$ between $(m_i)\in \ram$ and $(x_i)\in bb(e_i)\cup fbb(e_i)$:
$$
\ku R(m_i,x_i)\equi \a i\; \e j\;\;
m_0<x_i<m_j<m_{j+1}<x_{i+1}.
$$

The following is our basic determinacy result. The proof is based on
an idea already used in \cite{fr} to prove a different result and is
essentially descriptive set theoretical.
\begin{thm}\label{main}
Let $\A\subseteq bb(e_i)$ be analytic. Then the following are equivalent:
\begin{itemize}
\item[(a)] $\a \Delta>0\; \a (m_i)\;\e (x_i)\; (\ku R(m_i,x_i)\;\&\; (x_i)\in \A_\Delta)$,
\item[(b)] $\a \Delta>0$ II has a strategy in {\rm (IAG)} to play in $\A_\Delta$,
\item[(c)] $\a \Delta>0\; \e T$ block tree $([T]\subseteq \A_\Delta)$.
\end{itemize}
\end{thm}

\begin{proof}
The implication (c)$\saa$(a) is easy and (b)$\saa$(c) follows from Theorem
\ref{equivalence of games}, so we need only consider (a)$\saa$(b). Thus, assume
(a) and fix some $\Delta>0$. Let
$$
\B=\{(n_i)\in \ram\del
\e (x_i)\in \A_{\Delta/2}\; \a i\; n_{2i}<x_i<n_{2i+1}\}.
$$
Clearly $\B$ is analytic too. Moreover, we see that any $(m_i)\in \ram$
contains a subsequence in $\B$. For just choose some $(x_i)\in \A_{\Delta/2}$
such that  $\a i\; \e j\;\; m_0< x_i<m_j<m_{j+1}<x_{i+1}$. Then by leaving out
some terms of the sequence $(m_i)$ we get some subsequence $(n_i)$ such that
the indices match up. Using this fact, by Silver's theorem, there is some
infinite $B\subseteq \N$ such that $\ra{B}\subseteq \B$. And by the Jankov--von
Neumann selection theorem we can find a $\sigma({\bf \Sigma}^1_1)$-measurable
$f:\ra{B}\til bb(e_i)$ such that $f((n_i))=(x_i)$, where $(x_i)\in
\A_{\Delta/2}$ and $\a i\; n_{2i}<x_i<n_{2i+1}$.  Now choose inductively
sequences $(n_i)$, $(m_i)$ such that $n_i<m_i<n_{i+1}$ and sets $B_i\subseteq
B_{i-1}$ such that $m_{i-1}<\min B_i$ and such that for all $j_0<\ldots<j_k$
and
$$
C,D\in[(n_{j_0},m_{j_0},n_{j_1},m_{j_1},\ldots,n_{j_k},m_{j_k}),B_{j_k+1}]
$$
we have $\|f(C)_k-f(D)_k\|<{\delta_k}/2$. To see how this is done, start by
choosing $n_0<m_0$ in $B$ arbitrary and set $B_0=B$. Then for any $C\in
[(n_0,m_0),B_0]$ we have $f(C)_0\in \ku S({[e_{n_0+1},\ldots,e_{m_0-1}]})$,
and, as the unit sphere $\ku S({[e_{n_0+1},\ldots,e_{m_0-1}]})$ is compact, we
can by Silver's theorem find some $B_1\subseteq B_0$, $m_0<\min B_1$, such that
for all $C,D\in [(n_0,m_0),B_1]$ we have $\|f(C)_0-f(D)_0\|<{\delta_0}/2$. Now
suppose by induction that $B_l$ and $n_l<m_l$ in $B_l$ have been chosen. Then
we choose $B_{l+1}\subseteq B_l$ small enough that for all $j_0<\ldots<j_k\leq
l$ and all
$$
C,D\in
[(n_{j_0},m_{j_0},\ldots,n_{j_k},m_{j_k}),B_{l+1}]
$$
we have
$\|f(C)_k-f(D)_k\|<{\delta_k}/2$. Again this can be done as for all
$$
C,D\in [(n_{j_0},m_{j_0},\ldots,n_{j_k},m_{j_k}),B_{l}],
$$
we have
$$
f(C)_k,f(D)_k\in \ku S({[e_{n_{j_k}+1},\ldots,e_{m_{j_k}-1}]})
$$
and this unit sphere is compact. Moreover, we can
assume that $m_l<\min B_{l+1}$ and finally let $n_{l+1}<m_{l+1}$ be
arbitrary numbers in $B_{l+1}$.

It now remains to describe the strategy for II to play in $\A_\Delta$. First of
all, it is clear that we get an equivalent game if we demand that I plays a
subsequence $(n_{j_i})$ of the sequence $(n_i)$. Then II will respond to
$n_{j_0}<n_{j_1}<\ldots<n_{j_k}$ played by I with some $x_k$ satisfying
$n_{j_k}<x_k<m_{j_k}$ and such that for all $C\in
[(n_{j_0},m_{j_0},\ldots,n_{j_k},m_{j_k}),B_{j_k+1}]$ we have
$\|f(C)_k-x_k\|<{\delta_k}/2$. Thus, at the end of the game, when I has played
$(n_{j_i})$ and II has played $(x_k)$, we let $(y_i)=f((n_{j_i},m_{j_i})_i)$
and notice that for all $k$, $\|y_k-x_k\|<{\delta_k}/2$. As
$(y_i)=f((n_{j_i},m_{j_i})_i)\in \A_{\Delta/2}$, also $(x_i)\in\A_\Delta$. Thus
II will always play in $\A_\Delta$.\end{proof}

\begin{thm}\label{strat for I}
Let $\B\subseteq bb(e_i)$ be coanalytic. Then the following are equivalent:
\begin{itemize}
\item[(d)] $\e \Delta>0\; \e (m_i)\;\a (x_i)\; (\ku R(m_i,x_i)\til (x_i)\in {\rm Int}_\Delta(\B))$,
\item[(e)] $\e \Delta>0$ I has a strategy in {\rm (IAG)} to play in ${\rm Int}_\Delta(\B)$,
\item[(f)] $\e \Delta>0\; \a T$ block tree $([T]\cap {\rm
Int}_\Delta(\B)\neq\tom)$.
\end{itemize}
\end{thm}

\begin{proof}
Notice that if we define $\A$ in Theorem \ref{main} to be the set $\sim\!\B$,
then, as ${\rm Int}_\Delta(\B)=\sim(\A_\Delta)$,
$\neg$(d)$\equi$(a)$\equi$(c)$\equi\neg$(f). On the other hand,
(b)$\equi\neg$(e) only if the game is determined. But clearly,
(d)$\saa$(e)$\saa\neg$(b)$\saa\neg$(a)$\saa$(d). Thus (d)$\equi$(e)$\equi$(f).
\end{proof}

For the following result, we notice that the proof of Theorem
\ref{main} goes through for arbitrary projective sets as long as we
also assume a sufficient amount of determinacy. This follows from
the fact that projective determinacy implies that projective sets
are completely Ramsey and, moreover, that they can be uniformised by
projective sets (see \cite{kechris}).
\begin{thm}[Projective determinacy]
Let $\A\subseteq bb(e_i)$ be projective. Then the following are equivalent:
\begin{itemize}
\item[(g)] $\a \Delta>0\; \a (m_i)\;\e (x_i) \;(\ku R(m_i,x_i)\;\&\; (x_i)\in \A_\Delta)$,
\item[(h)] $\a \Delta>0$ II has a strategy in {\rm (IAG)} to play in $\A_\Delta$,
\item[(i)] $\a \Delta>0\; \e T$ block tree $([T]\subseteq \A_\Delta)$.
\end{itemize}
\end{thm}

With this result in hand, we can also get a perhaps more
satisfactory version of Theorem \ref{strat for I}.
\begin{thm}[Projective determinacy]\label{proj strat for I}
Let $\A\subseteq bb(e_i)$ be projective. Then the following are equivalent:
\begin{itemize}
\item[(j)] $\a \Delta>0\; \e (m_i)\;\a (x_i)\; (\ku R(m_i,x_i)\til (x_i)\in \A_\Delta)$,
\item[(k)] $\a\Delta>0$ I has a strategy in {\rm (IAG)} to play in $\A_\Delta$,
\item[(l)] $\a \Delta>0\; \a T$ block tree $([T]\cap \A_\Delta\neq\tom)$.
\end{itemize}
\end{thm}

\begin{proof}
Clearly, (j)$\saa$(k)$\saa$(l). Now suppose $\neg$(j) holds. Then
$$
\e \Delta_0>0\; \a (m_i)\;\e (x_i)\; (\ku R(m_i,x_i)\;\&\;
(x_i)\notin \A_{\Delta_0}).
$$
Now put $\B=\sim(\A_{\Delta_0})$, then a fortiori
$$
\a \Delta>0\; \a (m_i)\;\e (x_i)\; (\ku R(m_i,x_i)\;\&\; (x_i)\in
\B_{\Delta}),
$$
so also $\a \Delta>0\; \e T$ block tree $([T]\subseteq \B_\Delta)$. In
particular, if $\Delta={\Delta_0}/2$, then
$$
\e T \textrm{ block tree } ([T]\subseteq
\B_\Delta=(\sim\!(\A_{\Delta_0}))_{{\Delta_0}/2}\subseteq
\sim\!(\A_{{\Delta_0}/2})),
$$
 i.e., $\neg$(l), finishing the proof.
\end{proof}
There is a different and more constructive approach to Theorem \ref{main},
namely that of \cite{exact} which changes the position of the
$\Delta$-expansions. On the other hand, this method relies directly on the
determinacy of games and therefore only applies to Borel sets. We should also
mention that Odell and Schlumprecht \cite{os} prove a result concerning closed
sets that resembles Theorem \ref{proj strat for I}. Their exact result can also
be obtained by our methods.

\section{Applications}\label{applications}
In the following, we let $(e_i)$ be a Schauder basis for a Banach space. We
also restrict our attention to the infinite asymptotic game (IAG), which is
without consequence by Theorem \ref{equivalence of games}. Thus, I and II
always refer to the players in (IAG).

\begin{prop}\label{uniformity for II}
Suppose that $\A\subseteq bb(e_i)$ is analytic such that II has a strategy to
play in
$$
\A_\sim=\{(x_i)\in bb(e_i)\del \e (y_i)\in \A\; (x_i)\sim
(y_i)\}.
$$
Then there is some $K\geq 1$ such that II has a strategy to play in
$$
\A_K=\{(x_i)\in bb(e_i)\del \e (y_i)\in \A\; (x_i)\sim_K (y_i)\}.
$$
\end{prop}
\begin{proof}
To see this, we notice that if not, then, by Theorem \ref{main}, for
every $K$ there is some $(m_i^K)$ such that
$$
\a (x_i)\;(\ku R(m_i^K,x_i)\til (x_i)\notin \A_{c(K)\cdot K}),
$$
where $c(K)$ is a constant depending on the constant of the basis such that
whenever $(y_i)$ and $(x_i)$ are two elements of $bb(e_i)$ differing in at most
the $K$ first terms, then $(x_i)\sim_{c(K)}(y_i)$. Moreover, we can assume that
$m_i^K+1<m_{i+1}^K$. We use here that for $\Delta>0$ sufficiently small, we
have for all $K$, $(\A_K)_\Delta\subseteq \A_{K+1}$. Now pick an increasing
sequence $(k_i)$ such that
$$
\a K\; \a j\geq K\; \e l\geq 2j+1\; (k_j\leq m_l^K<m_{l+1}^K\leq
k_{j+1}).
$$
Suppose now that $(x_i)$ satisfies $\ku R(k_i,x_i)$. We claim that $(x_i)\notin
\A_\sim$. To see this, let $K$ be given. Then
$$
\a i\; \e j\geq i\; (x_i<k_j<k_{j+1}<x_{i+1}),
$$
so therefore we have also
$$
\a i\geq K\;\e l\geq 2i+1\; (x_i<m_l^K<m_{l+1}^K<x_{i+1}).
$$
Choose arbitrary $y_0,\ldots,y_K\in \ku S_E$ satisfying
$m_{2i}^K<y_i<m_{2i+1}^K$, and notice that then
$m_{2K}^K<y_K<m_{2K+1}^K<m^K_{2K+2}<x_{K+1}$. So the sequence
$$
(z_i)=(y_0,\ldots,y_K,x_{K+1},x_{K+2},x_{K+3},\ldots)
$$
belongs to $bb(e_i)$ and $(z_i)\sim_{c(K)}(x_i)$ as the two sequences only
differ by the first $K$ terms. On the other hand, we have $\ku R(m_i^K,z_i)$,
which implies that $(z_i)\notin \A_{c(K)\cdot K}$, and hence $(x_i)\notin
\A_K$. As $K$ is arbitrary, this show that $(x_i)\notin \A_\sim$ and thus that
$$
\a (x_i)\; (\ku R(k_i,x_i)\til (x_i)\notin \A_\sim).
$$
But then I has an obvious strategy to play in $\sim\!\A_\sim$, contradicting
that II should have a strategy to play in $\A_\sim$.
\end{proof}

The same argument easily shows

\begin{prop}
Suppose that $\A\subseteq bb(e_i)$ is analytic such that II has a strategy to
play in
$$
\A_{\iso}=\{(x_i)\in bb(e_i)\del \e (y_i)\in \A\; [x_i]\iso [y_i]\}.
$$
Then there is some $K\geq 1$ such that II has a strategy to play in
$$
\A_K=\{(x_i)\in bb(e_i)\del \e (y_i)\in \A\; [x_i]\iso_K [y_i]\}.
$$
\end{prop}

We now come to our first application of the results of Section \ref{evn}. Since
$(e_i)$ is a normalised Schauder basis of a Banach space $E$, we can consider
the dense subspace linearly spanned by the basis. When there is no chance of
confusion, we shall not insist of the difference between $E$ and this subspace.

\begin{prop}\label{rosenthal}
Let $(e_i)$ be a normalised Schauder basis for Banach space such that any block
tree has a branch equivalent to $(e_i)$. Then $(e_i)$ is equivalent to the
standard unit vector basis in either $c_0$ or $\ell_p$ for some $1\leq
p<\infty$.
\end{prop}

\begin{proof}
We claim that there is a normalised block basis $(x_i)\in bb(e_i)$ of $(e_i)$
that is perfectly homogeneous. For otherwise, we would for every $(x_i)\in
bb(e_i)$ have a further block basis $(y_i)\in bb(e_i)$ of $(x_i)$ such that
$(y_i)\not\sim(x_i)$. In particular, either $(x_i)\not\sim(e_i)$ or
$(y_i)\not\sim(e_i)$. But then it is easy to see that
$$
\a (m_i)\; \e (x_i)\; (\ku R(m_i,x_i)\;\&\;(x_i)\not\sim(e_i)),
$$
and hence by Theorem \ref{main}, there is a block tree $T$ all of whose
branches are inequivalent with $(e_i)$, contradicting the assumption.

So pick a perfectly homogeneous normalised $(x_i)\in bb(e_i)$, which by
Zippin's theorem (see \cite{ak}) is equivalent with either $c_0$ or some
$\ell_p$. On the other hand, by constructing a block tree of subsequences of
$(x_i)$, we see that also $(x_i)\sim(e_i)$.
\end{proof}

We recall  that a normalised basic sequence $(e_i)$ is called a {\em Rosenthal}
basic sequence if any normalised block has a subsequence equivalent to $(e_i)$
(see Ferenczi, Pelczar, and Rosendal \cite{fpr} for more on such bases). It is
still an open question whether Rosenthal sequences are equivalent to the
standard unit vector bases in $c_0$ or $\ell_p$, though the answer is positive
in case there is some uniformity or the subsequence can be chosen continuously.
The preceding proposition is in the same vein.

\begin{defi}
A {\em weakly null tree} is a non-empty set of finite strings of normalised
vectors $T\subseteq \ku S_E^{<\N}$ closed under initial segments such that for
all $(x_0,\ldots,x_{n-1})\in T$ the set
$$
\{y\in \ku S_E\del (x_0,\ldots,x_{n-1},y)\in T\}
$$
can be written as $\{y_i\}_{i=0}^\infty$ for some weakly null
sequence $(y_i)$.
\end{defi}

We next consider an analogous situation but for isomorphism instead
of equivalence.

\begin{thm}[Projective determinacy]\label{block homo}
Let $E$ be an infinite-dimensional Banach space that is not $\ell_1$-saturated. Then any weakly null
tree in $E$ has a branch spanning a space isomorphic to $E$ if and
only if $E$ has an unconditional block homogeneous basis, i.e., an
unconditional basis all of whose blocks span isomorphic spaces.
\end{thm}

\begin{proof}
Suppose first that any weakly null tree in $E$ has a branch spanning a space
isomorphic to $E$ and let $X$ be a separable infinite-dimensional subspace not
containing $\ell_1$. Then by Gowers' dichotomy theorem \cite{gowers}, $X$
either contains an unconditional basis or a hereditarily indecomposable (HI)
subspace generated by a basis $(e_i)$. In the latter case, by Rosenthal's
$\ell_1$-theorem, we can find a weakly Cauchy subsequence $(e'_i)$ and thus by
looking at the block sequence
$(f_i)_{i=0}^\infty=(e'_{2i}-e'_{2i+1})_{i=0}^\infty$, we get a weakly null
basic sequence. Let $T$ be the weakly null tree whose branches are exactly the
subsequences of $(f_i)_{i=0}^\infty$. Then $T$ has a branch spanning a space
isomorphic to $E$ and hence $E$ is isomorphic to a proper subspace that is HI,
contradicting the properties of HI spaces \cite{gowersmaurey}. Thus, $X$ must
contain an unconditional basic sequence $(g_i)$, which by James' theorem must
be shrinking, whereby any normalised block basis is weakly null. This means
that any block tree $T$ over the basis $(g_i)$ is in particular a weakly null
tree and hence has a branch spanning a space isomorphic to  $E$. So by Theorem
\ref{proj strat for I} we can find an increasing sequence $(m_i)$  such that
for all normalised blocks $(x_i)\in bb(g_i)$ of $(g_i)$, if $\ku R(m_i,x_i)$,
then $[x_i]$ is isomorphic to $E$. In particular, $(g_i)$ has a subsequence
that is block homogeneous and spans a space isomorphic to $E$.

Conversely, assume $E$ has a block homogeneous basis $(e_i)$ (we do
not need it to be unconditional). Any weakly null tree in $E$ has a
branch equivalent to a block basis of $E$ and thus spanning a space
isomorphic to $E$.
\end{proof}
It is plausible that one should be able to avoid the use of
projective determinacy in the above theorem. However, as projective
determinacy seems to be a part of the {\em right} axioms for set
theory one should not be too reluctant in using it.

It is unclear which spaces can have a block homogeneous basis. It
seems quite likely that this should only happen for the spaces
$\ell_p$ and $c_0$, but the question appears to be wide open. Modulo
the non-trivial fact that $\ell_p$, $p\neq 2$, and $c_0$ are not
homogeneous, a positive answer would of course provide another
solution to the homogeneous space problem.

\section{Subspaces of spaces with F.D.D.'s}\label{subspaces of
spaces with FDD's}

\subsection{Approximate games}
We would now like also to consider spaces that do
not necessarily have a basis or even an F.D.D. and thus obtain
coordinate free versions of Theorem \ref{main} and \ref{strat for
I}. However, when $E$ is a closed subspace of a space $F$ with an
F.D.D. $(F_i)$, we cannot immediately apply the results developed
till now. The problem is that there might be relatively few {\em
blocks} in $E$ for II to play and hence this would put unnatural
restrictions on the play of II. Instead we will need to consider
approximate blocks and thus continuously work with an extra layer of
approximations. This does not make for a nicer theory, and indeed
obscures many of the arguments that are otherwise very similar to
the previous setting. Unfortunately, this seems to be unavoidable
for certain applications.

We fix in the following spaces $E\subseteq F$ and an F.D.D. $(F_i)$
of $F$. We let for each interval $I\subseteq \N$, $P_I$ denote the
projection from $F$ onto the subspace $\oplus_{i\in
I}F_i=\overline{\rm span}(\bigcup_{i\in I}F_i)$.
\begin{defi}
We denote by $\Delta>0$ the fact that
$\Delta=(\delta_i)_{i=0}^\infty$ for some decreasing sequence of
$\delta_i>0$ converging to $0$. Given $\Delta>0$, a finite or
infinite sequence $(x_i)$ of vectors $x_i\in \ku S_E$ is said to be
a $\Delta$-{\em block} if there are intervals $I_i\subseteq \N$ such
that
$$
I_0<I_1<I_2<\ldots
$$
and for every $i$,
$$
\|P_{I_i}(x_i)-x_i\|<\delta_i.
$$
We also write $\ku B_\Delta(x_i,I_i)$ to denote that $(x_i)$ is a
$\Delta$-block as witnessed by $(I_i)$ and  denote by
$bb_{E,\Delta}(F_i)$ the set of $\Delta$-blocks (in $\ku
S_E^\infty$) with respect to the decomposition $(F_i)$ of $F$.
\end{defi}
We notice that due to the convergence of the $\delta_i$'s, if
$(x_i)$ is an infinite $\Delta$-block and $\Gamma>0$, then $(x_i)$
has an infinite subsequence which is a $\Gamma$-block. We should
also mention that the sequence $(I_i)$ witnessing that $(x_i)$ is a
$\Delta$-block is not necessarily unique, but if it exists then
there is a minimal one, i.e., such that $\max I_i\leq \max J_i$ for
any other witness $(J_i)$. This is easily seen by constructing the
$I_i$'s inductively.

Furthermore, if $K$ is the constant of the decomposition $(F_i)$ and
$(x_i)$ and $(y_i)$ are given such that $\a i\;
\|x_i-y_i\|<\delta_i$, then if $(x_i)$ is a $\Delta$-block, $(y_i)$
is a $4K\Delta$-block.  For if $(I_i)$ witnesses that $(x_i)$ is a
$\Delta$-block then for every $i$
\begin{displaymath}\begin{split}
\|y_i-P_{I_i}(y_i)\|&\leq
\|y_i-x_i\|+\|x_i-P_{I_i}(x_i)\|+\|P_{I_i}(x_i)-P_{I_i}(y_i)\| \\
&\leq \delta_i+\delta_i+2K\|x_i-y_i\| \\
&<4K\delta_i.
\end{split}\end{displaymath}

\begin{defi}
A $\Delta$-{\em block tree} $T$ is a non-empty subset $T\subseteq (\ku
S_E)^{<\N}$ closed under initial segments such that for all
$(x_0,\ldots,x_{n-1})\in T$ the set
$$
\{y\in \ku S_E\del (x_0,\ldots,x_{n-1},y)\in T\}
$$
can be written as $\{y_i\}_{i=0}^\infty$, where for each $i$ there
is an interval $I_i\subseteq \N$ satisfying
\begin{itemize}
  \item $\|y_i-P_{I_i}(y_i)\|<\delta_n$,
  \item $\min I_i\Lim{i\til \infty}\infty$.
\end{itemize}
\end{defi}
It is easily seen that any $\Delta$-block tree $T$ has a $\Delta$-spreading
subtree $S\subseteq T$ such that any branch $(x_i)\in [S]$ is a $\Delta$-block.

\begin{defi}
If $(m_i)\in \ram$ and $(I_i)$ are intervals, we write $\ku
R(m_i,I_i)$ in case
$$
m_0<I_0\;\&\; \a i\;\e j\; I_i<m_j<m_{j+1}<I_{i+1}.
$$
Similarly, if $(x_i)\in \ku S_E^\infty$ and $\Delta>0$, we write
$\ku R_\Delta(m_i,x_i)$ if $(x_i)$ is a $\Delta$-block with witness
$(I_i)$ satisfying $\ku R(m_i,I_i)$.
\end{defi}

Proceeding with our definitions, we let the {\em (infinite
asymptotic) $\Delta$-game} between two players I and II be defined
as follows: I and II alternate (with I beginning) in choosing resp.
natural numbers $n_0, n_1, n_2,\ldots$ and normalised vectors
$x_0,x_1,x_2,\ldots\in E$ and intervals $I_i$ such that for each
$i$, $\|x_i-P_{I_i}(x_i)\|<\delta_i$ and $n_i<I_i<I_{n+1}$.
$$
\begin{array}{cccccccccccc}
{\bf I} & & & n_0 &  & n_1 &  & n_2 & & n_3& &\ldots \\
{\bf II} & & &  & (x_0,I_0) &  & (x_1,I_1) &  & (x_2,I_2) & &
(x_3,I_3) &\ldots
\end{array}
$$

Finally, we now define for each set $\A\subseteq \ku S_E^\infty$ and
$\Delta>0$, the sets
$$
\A_\Delta=\{(y_i)\in\ku S_E^\infty\del \e (x_i)\in \A\; \a i\;
\|x_i-y_i\|<\delta_i\},
$$
and
\begin{displaymath}\begin{split}
{\rm Int}_\Delta(\A)=&\{(y_i)\in \ku S_E^\infty\del \a (x_i)\in \ku S_E^\infty\;( \a i\; \|x_i-y_i\|<\delta_i\til (x_i)\in \A)\}\\
=&\sim(\sim\A)_\Delta.
\end{split}\end{displaymath}
So as before ${\rm Int}_\Delta(\A)\subseteq \A\subseteq \A_\Delta$.
Moreover, there is some $\Delta>0$ depending only on the constant
$K$ of the decomposition $(F_i)$ such that if $\A\subseteq \ku
S_E^\infty$ is closed under taking equivalent sequences then
$\A_\Delta\cap bb_{E,\Delta}(F_i)\subseteq {\rm Int }_\Delta(\A)$.

\begin{thm}Suppose spaces $E\subseteq F$ are given, where $F$ has an
F.D.D $(F_i)$. Let $\A\subseteq \ku S_E^\infty$ be analytic. Then
the following are equivalent.
\begin{itemize}
  \item [(a)] $\a \Delta>0\;\a (m_i)\;\e (x_i)\; \big(\ku R_\Delta (m_i,x_i)\;\&\; (x_i)\in \A_\Delta\big)$.
  \item [(b)] $\a\Delta>0$ II has a strategy in the  $\Delta$-game to play into $\A_\Delta$.
  \item [(c)] $\a\Delta>0$ $\e T$ $\Delta$-block tree such that
  $[T]\subseteq \A_\Delta$.
\end{itemize}
\end{thm}

The proof goes along the lines of the proof of Theorem \ref{main}.

\begin{proof}
It should be obvious that from a $\Delta$-block tree all of whose branches
belong to $\A_\Delta$ one easily obtains a strategy for II in the $\Delta$-game
to play in $\A_\Delta$. Conversely, as in the proof of Proposition
\ref{equivalence of games}, any such strategy for II gives rise to a
$\Delta$-block tree with branches only in $\A_\Delta$. Thus (b)$\equi$(c) and
trivially (b)$\saa$(a).

For (a)$\saa$(b), assume that (a) holds and fix some $\Delta>0$.
Let
$$
\B=\{(n_i)\del \e(x_i)\;\e(I_i)\;\big( \ku
B_{\Delta/6}(x_i,I_i)\;\&\;\a i\; (n_{2i}<I_i<n_{2i+1})\;\&\;
(x_i)\in \A_{\Delta/6}\big)\}.
$$
Then as in the proof of Theorem \ref{main}, we find some infinite
$B\subseteq \N$ and $C$-measurable functions $f,g$ that to each
$(n_i)\in \ra{B}$ associate $(I_i)=g((n_i))$ and $(x_i)=f((n_i))$
such that $\ku B_{\Delta/6}(x_i,I_i)$, $\a i\;
(n_{2i}<I_i<n_{2i+1})$, and  $(x_i)\in \A_{\Delta/6}$.

Again construct inductively sequences $n_i<m_i<n_{i+1}$ and sets
$B_{k+1}\subseteq B_k$ such that $m_{i-1}<\min B_i$ and such that
for all $j_0<\ldots<j_k$ and all
$$
C,D\in [(n_{j_0},m_{j_0},\ldots,n_{j_k},m_{j_k}),B_{k+1}]
$$
we have $g(C)_k=g(D)_k$ and
$$
\|P_{g(C)_k}(f(C)_k)-P_{g(D)_k}(f(D)_k)\|<\delta_k/6.
$$
This can be done since in this case $n_{j_k}<g(C)_k<m_{j_k}$ and
there are only finitely many intervals $I$ such that
$n_{j_k}<I<m_{j_k}$, and as furthermore the ball of radius $2K$ in
$\oplus_{i\in g(C)_k}F_i$ is compact.

Finally, we describe the strategy for II to play in $\A_\Delta$
assuming that I is playing a subsequence of $(n_i)$. So suppose I
has played $n_{j_0}<\ldots<n_{j_k}$. Then II responds with
$(x_k,I_k)$ such that $\|P_{I_k}(x_k)-x_k\|<\delta_k/6$ and for all
$$
C\in [(n_{j_0},m_{j_0},\ldots,n_{j_k},m_{j_k}),B_{k+1}]
$$
we have $g(C)_k=I_k$ and
$$
\|P_{I_k}(f(C)_k)-P_{I_k}(x_k)\|<\delta_k/6.
$$
Thus, if at the end of the game, I has played $(n_{j_i})$ and II has
played $(x_i)$, then if $(y_i)=f((n_{j_i},m_{j_i})_i)$, we have for
all $i$
\begin{displaymath}\begin{split}
\|y_i-x_i\|&\leq
\|y_i-P_{I_i}(y_i)\|+\|P_{I_i}(y_i)-P_{I_i}(x_i)\|+\|P_{I_i}(x_i)-x_i\| \\
&\leq \delta_i/6+\delta_i/6+\delta_i/6 \\
&=\delta_i/2.
\end{split}\end{displaymath}
As $(y_i)=f((n_{j_i},m_{j_i})_i)\in \A_{\Delta/6}$, also $(x_i)\in
\A_\Delta$.
\end{proof}

\begin{cor}\label{strategy for I}
Let $\B\subseteq \ku S_E^\infty$ be coanalytic. Then the following
are equivalent:
\begin{itemize}
\item[(d)] $\e \Delta>0\; \e (m_i)\;\a (x_i)\; (\ku R_\Delta(m_i,x_i)\til (x_i)\in {\rm Int}_\Delta(\B))$,
\item[(e)] $\e \Delta>0$ I has a strategy in the $\Delta$-game to play in ${\rm Int}_\Delta(\B)$,
\item[(f)] $\e \Delta>0\; \a T$ $\Delta$-block tree $([T]\cap {\rm Int}_\Delta(\B)\neq\tom)$.
\end{itemize}
\end{cor}

\subsection{Uniformity}

Before we state the next proposition, which is very similar to a
result of Odell, Schlumprecht, and Zs\'ak \cite{osz}, let us recall
a few notions. If $(x_i)$ and $(y_i)$ are two basic sequences, we
say that $(y_i)$ {\em dominates} $(x_i)$ if for some $K$ and for all
choices of scalars $a_0,\ldots,a_n$ we have
$\|\sum_{i=0}^na_ix\|\leq K\|\sum_{i=0}^na_iy_i\|$. We denote this
by $(x_i)\leq(y_i)$ or $(x_i)\leq^K(y_i)$ if it holds for the
constant $K$.

\begin{prop}\label{uniformity for I}
Let $E$ be a closed infinite dimensional subspace of a Banach space $F$ with an
F.D.D. $(F_i)$, let $\Delta>0$, and suppose that there is a countable set $\A$
of normalised Schauder bases such that any $\Delta$-block tree has a branch
that is majorised by some element of $\A$. Then there is a constant $C$, a
sequence $(e_i)\in \A$, and a sequence $(m_i)$ such that whenever $(x_i)\in
bb_{E,\Delta}(F_i)$ satisfies $\ku R_\Delta(m_i,x_i)$ then $(x_i)\leq^C(e_i)$.

Similarly for $\geq $ and $\approx$.
\end{prop}

\begin{proof}
Notice first that by Corollary \ref{strategy for I}, we can find some sequence
$(m_i)$ such that
$$
\a (x_i)\in bb_{E,\Delta}(F_i)\;(\ku R_\Delta(m_i,x_i)\til \e
(f_i)\in \A\; (x_i)\leq (f_i)).
$$
By passing to a subsequence of $(m_i)$, we can suppose that for
every $i$ there is some $v\in \ku S_E$ such that
$\|v-P_{]m_i,m_{i+1}[}(v)\|<\delta_i$. List now $\A$ as
$\{(y_i^0),(y_i^1),(y_i^2),\ldots\}$ and suppose towards a
contradiction that the conclusion of the Proposition fails for the
specific sequence $(m_i)$. Then for every $n$ and $C$ there is some
$\Delta$-block $(x_i^{n,C})$ such that $\ku R_\Delta(m_i,x^{n,C}_i)$
and
$$
(x_i^{n,C})\not\leq^C(y^n_i).
$$
But then for all $n$, $C$, and $M\leq N$ there is some $L$ big
enough such that
$$
(x_i^{n,L})_{i=M+N}^\infty\not\leq^{C}(y^n_i)_{i=M+N}^\infty.
$$
Notice that, as $\ku R_\Delta(m_i,x^{n,L}_i)$, there is an interval
$I$ such that $m_{2M+2N}<I$ and
$\|x^{n,L}_{M+N}-P_I(x_{M+N}^{n,L})\|<\delta_{M+N}$. So  find
$v_M,\ldots,v_{M+N-1}\in \ku S_E$ and intervals
$I_M<\ldots<I_{M+N-1}$ such that
\begin{displaymath}\begin{split}
&m_N<I_M<m_{N+1}<m_{N+2}<I_{M+1}<m_{N+3}<\ldots\\
&<m_{2N-2}<I_{M+N-1}<m_{2N-1}<m_{2N}\leq m_{2M+2N}<I
\end{split}\end{displaymath}
and $\|v_i-P_{I_i}(v_i)\|<\delta_{i}$.

Let now
$$
z_i=\left\{\begin{array}{ccc}
v_{i} & {\rm for } & M\leq i<M+N\\
x_i^{n,L} & {\rm for } & M+N\leq i<\infty
\end{array}\right..
$$
Then $(z_i)_{i=M}^\infty\not\leq^{C}(y^n_i)_{i=M}^\infty$, while
$(z_i)_{i=M}^\infty$ is a $(\delta_i)_{i=M}^\infty $-block for some
witnessing sequence $(I_i)_{i=M}^\infty$ of intervals such that
$m_N<I_M$ and
$$
\a i\geq M\; \e j\; (I_i<m_j<m_{j+1}<I_{i+1}).
$$
Cutting $(z_i)_{i=M}^\infty$ off at some finite $P$ we get the
following claim.

\noindent {\em Claim:} For all $n$, $C$, and $M\leq N$ there is a
sequence $(z_i)_{i=M}^P$ and intervals $(I_i)_{i=M}^P$ such that
\begin{itemize}
  \item $m_N<I_M$,
  \item $\a i\in [M,P]\;\; \e j\;\;\; (I_i<m_j<m_{j+1}<I_{i+1})$,
  \item $\a i\in [M,P]\;\;\|z_i-P_{I_i}(z_i)\|<\delta_i$,
  \item $(z_i)_{i=M}^P\not\leq^{C}(y^n_i)_{i=M}^P$.
\end{itemize}

We can now construct a sequence $(p_i)$ and some $\Delta$-block
$(z_i)$ such that $\ku R_\Delta(m_i,z_i)$ and for all $l$
$$
(z_i)_{i=0}^{p_{l}-1}\not\leq^{l}(y^{\pi(l)}_i)_{i=0}^{p_{l}-1},
$$
where $\pi:\N\til \N$ is some surjection hitting each number
infinitely often. To begin, let $p_0=0$, $n=\pi(0)$, $C=0$, $M=p_0$,
and $N=M$ and find some $(z_i)_{i=p_0}^P$ and $(I_i)_{i=p_0}^P$
according to the claim. Set $p_1=P+1$.

Assume now that $p_0,\ldots, p_l$ and $z_0,\ldots,z_{p_l-1}$ have
been chosen. Let $n=\pi(l)$, $C=l$, $M=p_l$, and let $N\geq M$ be
large enough such that $I_{p_l-1}<m_{N-1}<m_N$. Then we can find
some $(z_i)_{i=p_l}^P$ and $(I_i)_{i=p_l}^P$ according to the claim
and finally let $p_{l+1}=P+1$. This finishes the inductive
construction.

We then see that $\ku R_\Delta(m_i,z_i)$, while on the other hand
$$
(z_i)_{i=0}^{\infty}\not\leq(y^n_i)_{i=0}^{\infty}
$$
for any $n$, which is a contradiction.
\end{proof}

One should contrast the proof of Proposition \ref{uniformity for I}
with that of Proposition \ref{uniformity for II}. In the former the
diagonalisation is over the sequence $(m_i)$ while in the latter the
diagonalisation is over the sequence $(x_i)$.

\subsection{Introduction of coordinates}
To introduce coordinates in a space $E$ that does not have an F.D.D we shall
use the setup of Odell and Schlumprecht from \cite{os}. They developed much of
the theory used to take care of spaces without F.D.D., but only proved the
basic determinacy result in the case of closed games.

The fundamental result that allows us to introduce coordinates is the following
lemma from \cite{os} based on a result of W. B. Johnson, H. Rosenthal and M.
Zippin \cite{jrz}.

\begin{lemme}\label{johnson}
Let $E$ be a separable Banach space and $(Y_n)_{n=1}^\infty$ a sequence of closed
subspaces each having finite codimension in $E$. Then $E$ is
isometrically embeddable into a space $F$ having an F.D.D. $(F_i)_{i=0}^\infty$
such that, when identifying $E$ with its image in $F$, the following
holds
\begin{enumerate}
  \item\label{i} $E\cap {\rm span}(\bigcup_{i=0}^\infty F_i)$ is dense in $E$.
  \item For every $n$, the finite codimensional subspace $E_n=E\cap\overline {\rm span}(\bigcup_{i=n}^\infty
  F_i)$ is contained in $Y_n$.
  \item There is a constant $c>1$ such that for every $n\geq 1$, there is
  a  finite set $D_n\subseteq \ku S(F^*_0\oplus\ldots\oplus F^*_{n-1})$
  such that for any $x\in E$ we have
  $$ \|x\|_{E/Y_n}=\inf_{y\in Y_n}\|x-y\|\leq c\max_{w^*\in
  D_n}w^*(x).$$
\end{enumerate}
From (\ref{i}) it follows that for each $n$, $E\cap {\rm
span}(\bigcup_{i=n}^\infty F_i)$ is dense in $E_n$.

Moreover, if $E$ has separable dual, $(F_i)$ can be chosen to be shrinking, and
if $E$ is reflexive, $F$ can also be chosen to be reflexive.
\end{lemme}

\noindent{\bf Remark:} We should mention that from the lemma above
it follows that for all $n$, $\delta>0$, and $x\in \ku
S_E\cap\big(\overline{\rm span}(\bigcup_{i=n}^\infty
F_i)\big)_\delta$, there is some $y\in \ku S_{Y_n}$ with
$\|x-y\|\leq 4\delta c$. The argument can be found on page 4095 in
\cite{os}.

\begin{defi}
Given a decreasing sequence of finite codimensional subspaces $(Y_n)$ of $E$ as
above, a $(Y_n)$-{\em block tree} is a non-empty subset $S\subseteq \ku
S_E^{<\infty}$ closed under initial segments such that for all
$(x_0,\ldots,x_{n-1})\in S$ the set
$$
\{y\in \ku S_E\del (x_0,\ldots,x_{n-1},y)\in S\}
$$
can be written as $\{y_i\}_{i=0}^\infty$, where for each $i$ there
is an $n_i$ satisfying
\begin{itemize}
  \item $y_i\in Y_{n_i}$,
  \item $n_i\Lim{i\til \infty}\infty$.
\end{itemize}
\end{defi}

We now sum up exactly the amount of knowledge we need from Lemma \ref{johnson}
in the following proposition.

\begin{prop}\label{wbjohnson}
Let $E$ be a separable reflexive Banach space. Then there is some reflexive
$F\supseteq E$ having an F.D.D. $(F_i)_{i=0}^\infty$  and a constant $c>1$ such
that whenever $\Delta>0$ and $T\subseteq \ku S_E^{<\N}$ is a $\Delta$-block
tree with respect to $(F_i)$, there is a weakly null tree $S\subseteq \ku
S_E^{<\N}$ such that
$$
[S]\subseteq [T]_{8\Delta c}\quad\&\quad[T]\subseteq [S]_{8\Delta c}.
$$
Moreover, $E$ contains a block sequence of $(F_i)$.
\end{prop}

\begin{proof}
Notice that as $E^*$ is separable, we can chose $Y_n=\bigcap_{i\leq n}{\rm
ker}\,\phi_i$ for a dense set $\{\phi_i\}$ of continuous functionals on $E$.
Then any $(Y_n)$-block tree is a  weakly null tree. Using these $Y_n$, we embed
$E$ isometrically into $F$ as in Lemma \ref{johnson}. Then we can talk about
$\Delta$-blocks etc. in $E$. By the remark above, we see that if $T$ is a
$\Delta$-block tree, then by replacing vectors one by one we can construct a
$(Y_n)$-block tree $S$ such that $[S]\subseteq [T]_{8\Delta c}$ and
$[T]\subseteq [S]_{8\Delta c}$.

The moreover part easily follows from (1) of Lemma \ref{johnson}.
\end{proof}

We are now ready for an application of the basic determinacy result,
which strengthens a result of Odell and Schlumprecht from \cite{os}.

\begin{thm}\label{os}
Let $E$ be a reflexive Banach space such that any weakly null tree
has a branch that is a subsymmetric (possibly conditional) basic
sequence. Then $E$ embeds into an $\ell_p$ sum of finite-dimensional
spaces for some $1<p<\infty$.
\end{thm}

\begin{proof} Use Proposition \ref{wbjohnson} to find a reflexive superspace $F\supseteq E$  as
described. Let $\B=\{(x_i)\in \ku S_E^\infty\del (x_i) \textrm{ is a
subsymmetric basis }\}$, which is coanalytic (in fact $F_\sigma$), and find
$\Delta>0$ small enough such that $\B_{8\Delta c}\cap bb_{E,\Delta
}(F_i)\subseteq {\rm Int }_{\Delta}(\B)$. Then if $T$ is a $\Delta$-block tree,
all of whose branches are $\Delta$-blocks, we can find a weakly null tree $S$
such that $[S]\subseteq [T]_{8\Delta c}$. Thus, as $S$ has a branch in $\B$,
$T$ has a branch in ${\rm Int}_\Delta(\B)$.

We now apply Corollary \ref{strategy for I} in order to find some
$\Gamma>0$ and a sequence $(m_i)$ such that for any $\Gamma$-block
$(x_i)$, if $\ku R_\Gamma(m_i,x_i)$, then $(x_i)$ is a subsymmetric
basic sequence.

We claim that there is some basic sequence $(e_i)$ such that if $\ku
R_\Gamma(m_i,x_i)$ then $(x_i)\sim(e_i)$. To see this, notice that if $(x_i)$
and $(y_i)$ are given such that $\ku R_\Gamma(m_i,x_i)$ and $\ku
R_\Gamma(m_i,y_i)$, then we can find subsequences $(z_{2i})$ of $(x_i)$ and
$(z_{2i+1})$ of $(y_i)$ such that $\ku R_\Gamma(m_i,z_i)$. But then each of
$(z_i)$, $(x_i)$, and $(y_i)$ is subsymmetric and hence
$$
(x_i)\sim (z_{2i})\sim(z_i)\sim(z_{2i+1})\sim(y_i).
$$
Thus $(x_i)$ and $(y_i)$ are equivalent to some common $(e_i)$. We claim that
$(e_i)$ is perfectly homogeneous. To see this, pick a block sequence $(y_i)$ of
the decomposition $(F_i)$ such that each term $y_i$ belongs to $\ku S_E$ and
such that $\ku R(m_i,y_i)$. Then any normalised block sequence $(z_i)$ of
$(y_i)$ also satisfies $\ku R(m_i,z_i)$ and hence also $\ku R_\Gamma(m_i,z_i)$.
So both $(y_i)$ and all of its normalised block sequences are equivalent with
$(e_i)$ and  the latter is therefore perfectly homogeneous. By Zippin's Theorem
and since $E$ and hence $[e_i]$ is reflexive, $(e_i)$ is equivalent with the
unit vector basis $(f_i)$ in some $\ell_p$, $1<p<\infty$.

By Proposition \ref{uniformity for I} we find that there must be a constant
$C>0$ such that whenever $(x_i)$ is a $\Gamma$-block sequence such that $\ku
R_\Gamma(m_i,x_i)$, then $(x_i)\sim_C(f_i)$. Letting
$H_i=\oplus_{m_{2i-1}+1}^{m_{2i+1}}F_i$, where $m_{-1}=-1$, we see then that
for any sequence $(u_i)$ in $\ku S_E$ such that for some sequence $(t_i)$
\begin{equation}\label{1}
\|P_{\oplus_{j=t_{i-1}+1}^{t_i-1}H_j}(u_i)-u_i\|<\gamma_i
\end{equation}
we have $(u_i)\sim_{C'}(f_i)$ for some constant $C'$. (There is a hidden use of
subsymmetry of $(f_i)$ used here in order not to worry about a shift in the
indices.)

We are now in a position to finish our proof as in the proof of Theorem 4.1 in \cite{os}.
\end{proof}

\section{Subsequence extraction}\label{subsequence extraction}
It is interesting to see that the fact that every block tree has a branch with
a certain property is equivalent to I having a strategy to choose subsequences
in real time. We propose here to weaken the conditions on I by instead letting
her hesitate for a while or even to change her mind a couple of times. This
leads to less restrictive notions of subsequence selection that might be of
interest elsewhere.

We begin by defining the two games to be studied. Fix a Banach space $E$ with a
basis $(e_i)$.

The {\em $0$'th subsequence game} ${\rm(SG)}_0$ is defined as follows. Player I
and II alternate in choosing respectively finite strings $s_1,s_2,\ldots$ such
that $s_i\in \{\#,0,1\}^i$ and normalised blocks $x_1<x_2<\ldots\in E$.
$$
\begin{array}{cccccccccccc}
{\bf I} & &  &  & s_1 &  & s_2 & & s_3 &\ldots\\
{\bf II} & &   & x_1 &  & x_2 & & x_3& &\ldots
\end{array}
$$
Moreover, we demand that the $s_i$ satisfy the following coherence
condition. For each $n$, the sequence $s_{n+1}(n)$, $s_{n+2}(n)$,
$s_{n+3}(n)$, etc. will begin by a finite number of $\#$'s and then
followed by only $0$'s or only $1$'s. Thus, if we interpret
$s_i(n)=\#$ as I not yet having decided whether $x_n$ should belong
to the subsequence she is choosing, we see that she is allowed to
hesitate for a finite time, but must then decide once and for all.

We therefore see II as constructing an infinite normalised block basis
$(x_i)_{i\in \N}$, while I chooses a subsequence $(x_i)_{i\in A}$ by
letting $n\in A\equi \lim_{i\til\infty}s_i(n)=1$.

It is not hard to see that any strategy naturally provides a {\em continuous}
function $\phi:bb(e_i)\til \ram$ choosing subsequences.

The {\em $1$'st subsequence game} ${\rm(SG)}_1$ is defined as similarly. Player
I and II alternate in choosing respectively finite strings $s_1,s_2,\ldots$
such that $s_i\in \{0,1\}^i$ and normalised blocks $x_1<x_2<\ldots\in E$.
$$
\begin{array}{cccccccccccc}
{\bf I} & &  &  & s_1 &  & s_2 & & s_3 &\ldots\\
{\bf II} & &   & x_1 &  & x_2 & & x_3& &\ldots
\end{array}
$$
Moreover, we demand that the $s_i$ only oscillate finitely, i.e.,
for every $n$, $\lim_{i\til\infty} s_i(n)$ exists.

Again we see II as constructing an infinite normalised block basis $(x_i)_{i\in
\N}$, while I chooses a subsequence $(x_i)_{i\in A}$ by letting $n\in A\equi
\lim_{i\til\infty}s_i(n)=1$. However, in this case a strategy for I only
provides a {\em Baire class} 1 function $\phi:bb(e_i)\til \ram$ choosing
subsequences.

\

\noindent{\bf Example:} Consider the space $\ell_1\oplus \ell_2$
with its usual basis $(e_i)$ having norm
$$
\|\sum_ia_ie_i\|=\|\sum_ia_{2i+1}e_{2i+1}\|_1+\|\sum_ia_{2i}e_{2i}\|_2.
$$
We let
$$
\A=\{(y_i)\in bb(e_i)\del (y_i)\sim\ell_1 \textrm{ or }
(y_i)\sim \ell_2\}
$$
and
$$
\A_K=\{(y_i)\in bb(e_i)\del (y_i)\sim_K\ell_1 \textrm{ or }
(y_i)\sim_K \ell_2\}.
$$
We show that player I has a strategy in ${\rm (SG)}_1$ to choose
subsequences in $\A$. On the other hand, there is no $K$ such that
every normalised block basis has a subsequence in $\A_K$, and, in
fact, even if we demand that II plays blocks exclusively in $\ku
S_{\ell_1}\cup \ku S_{\ell_2}$, I still has no strategy in $({\rm
SG})_1$ to play in some $\A_K$.  Also, I does not have a strategy in
${\rm (SG)}_0$ to choose subsequences in $\A$.

Thus, in particular, there is no equivalent version of Proposition
\ref{uniformity for I} in the game $({\rm SG})_1$, i.e., we cannot
in general get uniformity results from the mere existence of a
winning strategy for I in the $1$'st subsequence game..

First, any normalised block $x\in \ell_1\oplus\ell_2$ can be written uniquely
as $\lambda y_1+(1-\lambda)y_2$ for normalised block vectors $y_1\in \ell_1$,
$y_2\in \ell_2$, and a scalar $\lambda\in [0,1]$. Moreover, suppose $(x_i)$ is
a normalised block sequence with corresponding scalars $(\lambda_i)$. Then if
the $\lambda_i$ are bounded away from $0$ by some $\eps>0$, the sequence
$(x_i)$ will be equivalent with $\ell_1$. However, the smaller we need to take
$\eps>0$, the worse the constant of equivalence with $\ell_1$ becomes and this
explains why there is no uniform constant for equivalence. On the other hand,
if $\lambda_i\rightarrow 0$ sufficiently fast, then $(x_i)$ is equivalent with
$\ell_2$. So fix some $\eps_i>0$ converging to $0$ sufficiently fast that if
$\lambda_i<\eps_i$ for all $i$, then $(x_i)\sim\ell_2$.

We now describe the strategy for I in response to a sequence $(x_i)$
with corresponding scalars $(\lambda_i)$ played by II. I will never
first exclude some $x_i$ from the subsequence and then later on
include it, so therefore, each coordinate will change value $0,1$ at
most once and this will be from $1$ to $0$.

I chooses all vectors $x_1,x_2,\ldots$ until she meets the first
$\lambda_{i_1}<\eps_1$. In this case, she will eliminate all
previous $x_i$ and only stick with $x_{i_1}$. Now she will continue
by choosing all subsequent $x_i$ until she meets the first
coordinate $i_2>i_1$ at which $\lambda_{i_2}<\eps_2$. She will then
suppress all the $x_i$ between $x_{i_1}$ and $x_{i_2}$ and choose
$x_{i_2}$. Again she will choose all further $x_i$ until she meets
the first $i_3>i_2$ such that $\lambda_{i_3}<\eps_3$, chooses this
$x_{i_3}$ and suppresses all $x_i$ between $x_{i_2}$ and $x_{i_3}$,
etc.

Thus at the end of the game, either there is some $n$ such that all
but finitely many chosen $\lambda_i$ are greater than $\eps_n$ or I
has chosen a subsequence $(x_{i_l})$ such that
$\lambda_{i_l}<\eps_l$ for all $l$. In the first case, the
subsequence is equivalent to $\ell_1$ and in the latter case
equivalent to $\ell_2$. So I has a strategy in ${\rm (SG)}_1$ to
play in $\A$.

Now assume that II is only allowed to play blocks in $\ku S_{\ell_1}\cup\ku
S_{\ell_2}$, i.e., such that the corresponding $\lambda$ is either $0$ or $1$,
and suppose towards a contradiction that for some $K\geq 1$, I has a strategy
in $({\rm SG})_1$ to play in $\A_K$. Then there is some $N\geq 0$ such that I
has a strategy to choose subsequences all of whose terms, except at most $N$,
belong to $\ell_1$, or all whose terms except at most $N$ belong to $\ell_2$.
We fix such a strategy $\sigma$ for I. So $\sigma$ is a function assigning to
any normalised block sequence of vectors from $\ku S_{\ell_1}\cup\ku
S_{\ell_2}$ a binary sequence of equal length.

We let $\PP$ be the set of finite block sequences $(x_1,\ldots,x_n)$
such that each term belongs to $\ku S_{\ell_1}\cup\ku S_{\ell_2}$.
Set $(x_1,\ldots, x_n)\geq (y_1,\ldots,y_m)$ if $n\leq m$ and
$x_1=y_1,\ldots,x_n=y_n$.

We  claim that for each $M$ and $p=1,2$ the set
$$
\E^p_M=\{(x_1,\ldots, x_j)\in \PP\del j> M\;\&\;\#(i\leq
M\;|\;\sigma(x_1,\ldots,x_j)(i)=1\;\&\; x_i\in \ell_p)\leq N\}
$$
is dense in $(\PP,\leq)$, i.e., any element of $\PP$ has a minorant in
$\E^p_M$. To see this, suppose $(x_1,\ldots,x_n)$, $p$ and $M$ are given.
Choose normalised consecutive blocks $x_{n+1},x_{n+2},\ldots$ belonging to
$\ell_{3-p}$ and let II play the infinite sequence $x_1,x_2,\ldots$. Since
$\sigma$ is a strategy in $({\rm SG})_1$, we know that $s=\lim_{j\til
\infty}\big(\sigma(x_1,\ldots,x_j)|_M\big)$ exists in $2^M$ and by assumption
on $\sigma$, we must have that
$$
\#(i\leq M\;|\;s(i)=1\;\&\; x_i\in \ell_p)\leq N.
$$
So fix some $j\geq n+M$ such that $\sigma(x_1,\ldots,x_j)|_M=s$,
then $(x_1,\ldots,x_j)\in \E^p_M$ and
$(x_1,\ldots,x_n)\geq(x_1,\ldots,x_j)$.

Let now $(x_1,x_2,x_3,\ldots)$ be an $\{\E_M^p\}_{p=1,2, M\in \N}$-generic
sequence, i.e., for all $M$ and $p$, $(x_i)$ has an initial segment belonging
to $\E^p_M$. Let also $A\subseteq \N$ be the ultimate response of II to this
sequence, i.e., $i\in A\equi \lim_{j\til \infty}\sigma(x_1,\ldots,x_j)(i)=1$.
Fix $p$ such that infinitely many terms in $(x_i)_{i\in A}$ belongs to $\ell_p$
and find $M$ large enough that there are $N+1$ many $i\leq M$ that belong to
$A$ and at the same time $x_i\in \ell_p$. Find $l>M$ such that for all $j\geq
l$, $\sigma(x_1,\ldots,x_j)|_M=\chi_A|_M$. But then if $j \geq l$ is such that
$(x_1,\ldots,x_j)\in \E^p_M$ we get a contradiction.

To see that I has no strategy in the $0$'th subsequence game to
choose subsequences in $\A$ is easier. For this it is enough to
notice that we can let II play vectors from $\ell_1$ until I commits
to take  at least one of these, then we let II continue to play
vectors from $\ell_2$ until I commits to at least one of these, and
II continues again with vectors from $\ell_1$ etc. Thus in the end,
I will have chosen infinitely many vectors from both $\ell_1$ and
$\ell_2$, so fails to play in $\A$.

\

It is proved in \cite{fpr} that if $(e_i)$ is a normalised basic
sequence such that every normalised block basis of $(e_i)$ has a
subsequence equivalent to $(e_i)$ (i.e., $(e_i)$ is a Rosenthal
basis), and, moreover, this subsequence can be chosen continuously
in the block basis, then $(e_i)$ is equivalent to the standard unit
vector basis of $c_0$ or some $\ell_p$. Thus, in particular, if I
has a strategy in $({\rm SG})_0$ to choose subsequences equivalent
to $(e_i)$, then $(e_i)$ is equivalent to $c_0$ or $\ell_p$. This
provides another proof of Proposition \ref{rosenthal}. A natural
question is whether this condition can be weakened to I having a
strategy in $({\rm SG})_1$.

\

\noindent{\em Acknowledgement:} I would like to thank Ted Odell,
Thomas Schlumprecht, and especially Valentin Ferenczi for many
helpful discussions on the topic of this paper.

\

\

\begin{flushleft}
{\em Address of C. Rosendal:}\\
Department of Mathematics, Statistics, and Computer Science\\
University of Illinois at Chicago\\
322 Science and Engineering Offices (M/C 249)\\
851 S. Morgan Street\\
Chicago, IL 60607-7045, USA\\
\texttt{rosendal@math.uic.edu}
\end{flushleft}

\end{document}